\documentclass[11pt]{amsart}
\newtheorem{thm}{Theorem}[section]
\newtheorem{lemma}[thm]{Lemma}

\newtheorem{conj}[thm]{Conjecture}
\newtheorem{prop}[thm]{Proposition}
\newcommand{\R}{\mathbb{R}}

\newcommand{\C}{\mathbb{C}}
\newcommand{\Z}{\mathbb{Z}}
\newcommand{\N}{\mathbb{N}}
\newcommand{\T}{\mathbb{T}}
\renewcommand{\H}{C^\alpha(\T)}
\newcommand{\q}{q}
\usepackage{amsmath}
\usepackage{graphicx}
\usepackage{hyperref}
\hypersetup{colorlinks=true, urlcolor=black, linkcolor=black, citecolor=black}

\begin{document}
\title[transfer operators on the circle with trigonometric weights]
{On transfer operators on the circle with trigonometric weights}
\date{\today}
\author{Xianghong Chen \and Hans Volkmer}
\address{X. Chen\\Department of Mathematical Sciences\\University of Wisconsin-Milwaukee\\Milwaukee, WI 53211, USA}
\email{chen242@uwm.edu}
\address{H. Volkmer\\Department of Mathematical Sciences\\University of Wisconsin-Milwaukee\\Milwaukee, WI 53211, USA}
\email{volkmer@uwm.edu}

\subjclass[2010]{37C30, 37E10, 35P15, 37N99}
\keywords{transfer operator, spectral problem, Cantor set}

\begin{abstract}
We study spectral properties of the transfer operators $L$ defined on the circle $\T=\mathbb R/\mathbb Z$ by
$$(Lu)(t)=\frac1 d\sum_{i=0}^{d-1} f\left(\frac{t+i}{d}\right)u\left(\frac{t+i}{d}\right),\ t\in\T$$
where $u$ is a function on $\mathbb T$. We focus in particular on the cases $f(t)=|\cos(\pi t)|^\q$ and $f(t)=|\sin(\pi t)|^\q$, which are closely related to some classical Fourier-analytic questions. We also obtain some explicit computations, particularly in the case $d=2$. Our study extends work of Strichartz \cite{Strichartz1990} and Fan and Lau \cite{FanLau1998}.
\end{abstract}
\maketitle

\section{Introduction}
Let $\T=\R/\Z$ and identify it with $[0,1]$ in the usual way. For integers $d\ge 2$, let $F$ be the $d$-adic Bernoulli map 
$$F: \mathbb T\rightarrow \mathbb T,\ t\mapsto d\cdot t\mod 1.$$
Let $f$ be a function (weight) on $\T$. Consider the weighted transfer operator $L$ associated to $F$:
\begin{align*}
(Lu)(t)
&=\frac1 d\sum_{s\in F^{-1}(t)} f(s)u(s)\\
&=\frac1 d\sum_{i=0}^{d-1} f\left(\frac{t+i}{d}\right)u\left(\frac{t+i}{d}\right),\ t\in\T
\end{align*}
where $u$ is a function on $\T$. Such operators are also called Ruelle (or Ruelle-Perron-Frobenius) operators. They can also be defined associated to more general maps and on more general spaces (cf. Hennion \cite{Hennion1993} and Baladi \cite{Baladi2000} for more background).

In this paper, we study spectral properties of $L$ as an operator acting on $C(\T)$, the space of continuous functions on $\T$. When $f\equiv 1$, this question has been extensively studied,
especially in the case $d=2$ (cf. Vep\v{s}tas \cite{V} and references therein). For more general weights $f$, there are Perron-Frobenius type theorems that describe spectral properties of $L$ (cf. \cite[Theorem~1.5]{Baladi2000}). However, such theorems often require $f$ to be strictly positive, which is not met by the main examples we are interested in:
$$(c)\ f(t)=|\cos(\pi t)|^\q \hspace{1cm} (s)\ f(t)=|\sin(\pi t)|^\q$$
where $\q>0$. In Section \ref{sec:quasicompact}, we develop Perron-Frobenius type theorems for transfer operators $L$ with such `degenerate' weights (more precisely, weights that have exactly one zero on $\T$). The theorems are derived using notions of quasicompactness and Krein property, which we verify by exploiting the specific structure of the Bernoulli map $F$; see also Fan and Lau \cite{FanLau1998} for similar treatments. As a corollary, we conclude that the operator $L$ satisfies classical Perron-Frobenius theorems in all cases of $d$ and $(c)/(s)$, except for the case $d=2$ and $(c)$. 

In Section \ref{sec:special-weights}, we study in more detail the spectral properties of $L$ in the non-exceptional cases. When $\q$ is an even integer, we obtain explicit computations of $\rho(L)$ (the spectral radius of $L$) by reducing to a finite-dimensional problem. When $\q$ is not an even integer, evaluating $\rho(L)$ is more difficult. We derive in this case estimations of $\rho(L)$, particularly for $d=3$ (note that when $d$ is odd, $(c)$ and $(s)$ are equivalent). As an application, we obtain asymptotic behavior of some integrals of the form
$$I_n=\int_{\T} \prod_{j=0}^{n-1} f(d^j t)dt,\quad \text{as } n\rightarrow\infty.$$
In particular, we extend a result of Strichartz \cite{Strichartz1990} concerning the Fourier transform of the middle-third Cantor set. We also study geometric properties of the function $L1$, 
as well as asymptotic behavior of $\rho(L)$ as $\q\rightarrow\infty$. For the latter question it turns out that one needs to distinguish the case when $d$ is even and $f$ is given by $(s)$.

In Section \ref{sec:binary}, we give a detailed account of the exceptional case $d=2$ and $(c)$. Using an explicit formula for the iterates $L^n 1$, we find the spectral radius and eigenfunctions of $L$ explicitly (see also Fan and Lau \cite{FanLau1998} for related results), and obtain geometric properties of $L^n1$ for $n\ge 1$ (especially for $\q\le 1$ and even $\q$'s). The spectral problem in this case is closely related to the case $f\equiv 1$ mentioned above, and has to do with the Hurwitz zeta functions.

In Section \ref{sec:Lp}, we study the spectral problem on Lebesgue spaces. In Section \ref{sec:application}, we give an application to Fourier multipliers.

\section{Quasicompact transfer operators}\label{sec:quasicompact}

Let $f:\R\to\R$ be a continuous nonnegative 1-periodic function, and let $d\ge 2$ be an integer.
We consider the transfer operator
\begin{equation}\label{L:L}
 (Lu)(t)=\frac1 d\sum_{i=0}^{d-1} f\left(\frac{t+i}{d}\right)u\left(\frac{t+i}{d}\right) .
\end{equation}
Let $\T=[0,1]$ with $0$ and $1$ identified (circle). Let $C(\T)$ be the Banach space of continuous complex-valued functions on $\T$
endowed with the maximum norm $\|\cdot\|_\infty$.
Then $L:C(\T)\to C(\T)$ is a bounded linear operator. Moreover, $L$ is positive in the sense that $u\ge 0$ implies $Lu\ge 0$.

Define a map $F:\T\to \T$ by $F(t)=d\cdot t \mod 1$. Then we can write
\[ (Lu)(t)=\frac1d \sum_{s\in F^{-1}(t)} f(s)u(s), \quad t\in\T .\]
For each $n\in\N$, set
\[ f_n(t)=\prod_{j=0}^{n-1} f\big(F^j(t)\big) .\]
Then
\begin{equation}\label{L:Ln}
(L^n u)(t)=\frac{1}{d^n} \sum_{s\in F^{-n}(t)} f_n(s)u(s) .
\end{equation}

Let $0<\alpha\le 1$. 
Consider the Banach space $\H$ 
of H\"older continuous functions $u:\T\to \C$ with the norm
\[ \|u\|_\alpha=\sup_{s\neq t} \frac{|u(s)-u(t)|}{|s-t|^\alpha} + \|u\|_\infty .\]
If $f\in \H$ then $L_\alpha u=Lu$ defines a bounded linear operator $L_\alpha:\H\to \H$.

Let $T$ be a bounded linear operator on a Banach space $X$. We denote its spectral radius by $\rho(T)$.
$T$ is called quasicompact if there exists a compact operator $K$ on $X$ such that $\rho(T-K)<\rho(T)$.
If $T$ is quasicompact and $\lambda\in\C$ is in the spectrum of $T$ with $|\lambda|>\rho(T-K)$, then $\lambda$ is an eigenvalue of
$T$.

The following theorem is proved in \cite[pages 3--4]{Sm}.  In \cite[Proposition 1]{Sm} it is assumed that $f$ is positive while we assume here that $f$ is nonnegative. However positivity of $f$ is not used on pages 3--4 of \cite{Sm}. See also \cite{Hennion1993}.

\begin{thm}\label{L:t1}
Let $0\le f\in \H$ for some $0<\alpha\le 1$.
Then $\rho(L)=\rho(L_\alpha)$. Furthermore, if  $\rho(L)>0$, then $L_\alpha$ is quasicompact.
\end{thm}

For each $n\in\N$, set
\begin{equation}\label{L:hn}
 h_n=L^n 1.
 \end{equation}
Define
\begin{equation}\label{L:rnRn}
 r_n=\min_{t\in\T}  h_n(t),\quad R_n=\max_{t\in\T} h_n(t) .
\end{equation}
It is easy to show that
\[ r_{m+n}\ge r_mr_n,\quad R_{m+n}\le R_mR_n .\]
Therefore, the limits
\begin{equation}\label{L:rR}
 r=\lim_{n\to\infty} (r_n)^{1/n} =\sup_n\ (r_n)^{1/n},\quad R=\lim_{n\to\infty} (R_n)^{1/n}=\inf_n\ (R_n)^{1/n}
\end{equation}
exist. In particular, for every $n\in\N$ we have
\begin{equation}\label{L:est1}
 r_n^{1/n}\le r\le R\le R_n^{1/n}.
\end{equation}
Moreover, since $R_n=\|L^n\|_{C(\T)\rightarrow C(\T)}$, by Gelfand's formula, $R=\rho(L)$.

\begin{thm}\label{L:t2}
Let $w\in C(\T)$ be a unit, that is, $w(t)>0$ for all $t\in\T$.
Then
\begin{equation}\label{L:est2}
\min_{t\in\T} \frac{(Lw)(t)}{w(t)}\le r\le R\le \max_{t\in\T} \frac{(Lw)(t)}{w(t)} .
\end{equation}
\end{thm}
\begin{proof}
We define a bounded linear operator $S$ on $C(\T)$ by
\[ (Su)(t)= \frac{1}{w(t)} L(wu)(t),\]
a sequence of functions
\[ \tilde h_n= S^n 1,\]
and sequences of numbers
\[ \tilde r_n =\min_{t\in\T} \tilde h_n(t),\quad \tilde R_n=\max_{t\in\T} \tilde h_n(t) .\]
Since $S$ is positive, we obtain for every $n\in\N$
\[
 (\tilde r_n)^{1/n}\le \lim_{n\to\infty}  (\tilde r_n)^{1/n}=: \tilde r\le \tilde R:= \lim_{n\to\infty}  (\tilde R_n)^{1/n}\le (\tilde R_n)^{1/n} .
\]
Since $w$ is a unite, there are constants $a,b>0$ such that
\[ a\le w(t)\le b\]
for all $t\in\T$. This implies
\[ a h_n(t)=(L^na)(t) \le w(t) \tilde h_n(t) \le (L^nb)(t)=b h_n(t) .\]
From this we obtain
\[  r_n\le \frac{b}{a} \tilde r_n,\quad \tilde r_n\le \frac{b}{a} r_n,\quad  R_n\le \frac{b}{a} \tilde R_n,\quad \tilde R_n\le \frac{b}{a} R_n .\]
Thus
\[ r=\tilde r,\quad R=\tilde R .\]
Now \eqref{L:est2} follows from
\[ \tilde r_1\le \tilde r=r\le \tilde R=R\le \tilde R_1.\]
\end{proof}

We say that $L$ is a Krein operator if, for all $u\in C(\T)$ such that $u(t)\ge 0$ for all $t\in\T$ but $u(t_0)>0$ for at least one $t_0\in\T$, there is $n\in\N$ such that $L^nu$ is a unit. Note that $n$ may depend on $u$.
It is easy to show that a Krein operator carries units to units (cf. \cite[Lemma 5.2]{A}).
It follows from \eqref{L:Ln} that if $f(t)>0$ for all $t$ then $L$ is a Krein operator,
Also, if $f$ vanishes on an interval of positive length, then $L$ cannot be a Krein operator.

\begin{lemma}\label{L:Krein1}
Suppose $f$ has exactly one zero in $[0,1)$.
If $f_n$ has four zeros that form an arithmetic progression with step size $d^{-n}$, then $d=2$ and $f(\frac12)=0$.
\end{lemma}
\begin{proof}
Let $s_0+\Z$ be the set of zeros of $f$.
Suppose that $t_i=t+id^{-n}$ is a zero of $f_n$ for $i=0,1,2,3$.
Then there exist integers $0\le k_i\le n-1$, and integers $j_i$ such that
\begin{equation}\label{eq}
 t_i=d^{-k_i}(s_0+j_i),\quad i=0,1,2,3.
\end{equation}
We will assume that $k_0=\max\{k_0,k_1,k_2,k_3\}$ (the other cases are mentioned at the end of the proof.)
Clearly, $k_0>k_1$.
Since we can replace $s_0$ by $s_0+j$ with any integer $j$, we will assume that $j_1=0$ in order to simplify the notation.
From $t_1-t_0=t_2-t_1=d^{-n}$, we obtain
\[ d^{-k_1}s_0-d^{-k_0}(s_0+j_0)=d^{-n},\quad d^{-k_2}(s_0+j_2) -d^{-k_1}s_0=d^{-n} .\]
Eliminating $s_0$ from these equations, we find
\[ (1-j_0d^{n-k_0})(d^{k_0-k_1}-d^{k_0-k_2})= (1-j_2d^{n-k_2})(1-d^{k_0-k_1}) .\]
This is an equation involving only integers. Since $n>k_2$, $k_0>k_1$, the right-hand side is not divisible by $d$.
Therefore, we must have $k_0=k_2$. But this is impossible when $d>2$. So we must have $d=2$ and $k_0=k_2=n-1$.

Now suppose that $d=2$ and $k_0=k_2=n-1$.
Without loss of generality we take $j_0=0$, $j_2=1$.
Since $t_1-t_0=t_3-t_2=2^{-n}$ we obtain
\[ 2^{n-k_1}(s_0+j_1)=2s_0+1,\quad 2^{n-k_3}(s_0+j_3)=2s_0+3 .\]
Eliminating $s_0$ this gives
\[ (3-j_32^{n-k_3})(2^{n-k_1-1}-1)=(1-j_12^{n-k_1})(2^{n-k_3-1}-1) .\]
It is clear that $k_1<n-1$, $k_3<n-1$. Therefore, the equation yields that
$2^{n-k_3-2}-3\cdot 2^{n-k_1-2}+1$ is an even integer. This implies that $k_1=n-2$ or $k_3=n-2$.
If $k_1=n-2$ then $2s_0+1=4(s_0+j_1)$, which implies that $s_0-\frac12$ is an integer.
If $k_3=n-2$, we have the same conclusion.

If $k_2=\max\{k_0,k_1,k_2,k_3\}$, the proof is almost the same. We again obtain $k_0=k_2=n-1$ and $d=2$ after the first part of the proof.
The other two cases can be reduced to the treated ones by replacing $s_0$ by $-s_0$.
\end{proof}

\begin{lemma}\label{L:Krein2}
Suppose $f$ has exactly one zero $s_0\in[0,1)$.
Then $L$ is a Krein operator unless $s_0=\frac12$ and $d=2$.
\end{lemma}
\begin{proof}
Let $u\in C(\T)$ be nonnegative but not identically zero,
Choose $n$ so large that $u(t)>0$ for $m d^{-n}\le t\le (m+4) d^{-n}$ for some integer $m$, $0\le m\le d^n-4$.
We claim that $(L^n u)(t)>0$ for all $t\in\T$. In fact, by Lemma \ref{L:Krein1}, if $t\in\T$ then among the four points $t_i=d^{-n}(t+m+i)$, $i=0,1,2,3$,
at least one satisfies $f_n(t_i)>0$. Then \eqref{L:Ln} implies $(L^n u)(t)\ge f_n(t_i)u(t_i)>0$.
\end{proof}

If $d=2$ and $f(\frac12)=0$ then $L$ is not a Krein operator: If $u(0)=0$ then $(L^nu)(0)=0$ for all $n\in\N$.

\begin{thm}[See also \cite{FanLau1998}]\label{L:t3}
Suppose that $L$ is a Krein operator. Then the following statements hold.\\
(a) $R>0$.\\
(b) If $Lv=\lambda v$ with $v\ne 0$ and $|\lambda|=R$,
then there is a constant $\theta\in\R$ such that $e^{-i\theta} v$ is a unit. \\
(c) $L$ has no eigenvalue $\lambda$ on the circle $|\lambda|=R$ except possibly $\lambda=R$.\\
(d) If $R$ is an eigenvalue of $L$, then its algebraic multiplicity is $1$.\\
(e) If $0\le f\in \H$, then  $L_\alpha$ is quasicompact, $r=R$, and $R$ is an eigenvalue of $L_\alpha$ of algebraic multiplicity $1$ with a unit eigenfunction. \end{thm}
\begin{proof}
(a) Since $L$ is a Krein operator, $h_1=L1$ is a unit, so $0<r_1\le R$.

\noindent (b)
$L v =\lambda v$ implies $z:=L|v|-R|v|\ge 0$.
Suppose that $z$ is not identically zero. Then there is $n\in\N$ such that $L^nz$ and $w:=L^n |v|$ are units.
It follows that there is $\delta>0$ such that
\[ (L^nz)(t)=(Lw)(t)-R w(t)\ge\delta w(t),\quad t\in \T.\]
Applying Theorem \ref{L:t2}, we obtain the contradiction
\[ R+\delta\le \min_{t\in\T} \frac{(Lw)(t)}{w(t)}\le R .\]
Therefore, $z=0$ and
\begin{equation}\label{L:eq}
L|v|=R|v|.
\end{equation}
Then
\[ |L^n v|=R^n |v|=L^n|v|=w.\]
Therefore, $|v|$ is a unit.
We claim that there is a constant $\theta\in\R$ such that $e^{-i\theta} v(t)>0$ for all $t$.
Suppose this is not true. Since $L|v|=|Lv|$, there is $n\in\N$ and $1\le i<j\le d^n$ such that $v(s)/v(t)\not\in(0,\infty)$ for all $s\in I_{n,i}$, $t\in I_{n,j}$.
Since $L$ is a Krein operator, $f$ does no vanish on an interval of positive length. Then also $f_n$ does not vanish on an interval of positive length.
Therefore, there is $s\in I_{n,i}$, $t\in I_{n,j}$ with $(t-s)d^{-n}\in\Z$ such that $f_n(s)\ne 0$, $f_n(t)\ne 0$.
Hence, by \eqref{L:Ln}, $|(L^n v)(s)|<(L^n |v|)(s)$, which is a contradiction.
Therefore, the claim is proved.

\noindent (c) follows from (b).

\noindent (d)
Suppose that $R$ is an eigenvalue of $L$. By (b), each corresponding eigenfunction is a constant multiple of a unit. It follows that
the geometric multiplicity of the eigenvalue $R$ is $1$.
Now assume that there are $u,w\in C(\T)$ such that $L u-Ru = w$, $L w=Rw$, where $w$ is a unit. We may assume that $u$ is a unit. There is $\delta>0$ such that $\delta u\le w$.
Then Theorem \ref{L:t2} leads to the contradiction
\[ R+\delta\le \min_{t\in\T} \frac{(Lu)(t)}{u(t)}\le R .\]
This shows that the algebraic multiplicity of the eigenvalue $R$ is $1$.

\noindent (e)
Since $\rho(L)=R>0$, it follows from Theorem \ref{L:t1} that $L_\alpha$ is quasicompact. Then $L_\alpha$ has an eigenvalue $\lambda$ on the circle $|\lambda|=R$.
By (c), $R$ is an eigenvalue of $L_\alpha$. There is a corresponding unit eigenfunction.
Now $r=R$ follows from Theorem \ref{L:t2}.
\end{proof}

\begin{thm}[See also \cite{FanLau1998}]\label{L:t4}
Suppose that $0\le f\in \H$ and that $L$ is a Krein operator.
Let $P$ be the spectral projection onto the eigenspace of $L_\alpha$ corresponding to the eigenvalue $R$.\\
(a) The sequence $R^{-n}L_\alpha^n$ converges to $P$ as $n\to\infty$ with respect to the operator norm.\\
(b) The sequence $R^{-n} h_n$ converges in $\H$ to an eigenfunction of $L_\alpha$ corresponding to the eigenvalue $R$.
\end{thm}
\begin{proof}
(a)
By Theorem \ref{L:t3}, $L_\alpha$ is quasicompact and the eigenvalue $R$ of $L_\alpha$ is an isolated point of its spectrum.
Therefore, there exists the spectral projection $P$ onto the one-dimensional root subspace belonging to the eigenvalue $R$.
The Banach space $\H$ is a direct sum of the subspaces $P\H$ and $(1-P)\H$. Both subspaces are invariant under $L_\alpha$.
On $P\H$, $L_\alpha$ acts as $R$ times the identity.
Set $S:=R^{-1}(1-P)L_\alpha$. By Theorem \ref{L:t3}, the spectral radius of $S$ is less than $1$ so $S^n$ converges to $0$ as $n\to\infty$ in the operator norm.
We have $R^{-n} L_\alpha^n =S^n+ P$ which implies statement (a).

\noindent (b)
We have $R^{-n}h_n=R^{-n}L^n_\alpha1\to P1$ as $n\to\infty$. Since
\[1\le R^{-n}R_n\le \|R^{-n} h_n\|_\alpha\to \|P1\|_\alpha,\]
we have $P1\ne0$.
\end{proof}

We consider now the following problem that was the original motivation for this paper.
Let $f:\R\to\C$ be a bounded measurable and 1-periodic function, and let $d\ge 2$ be an integer. For $n\in\N$, define as before
\[ f_n(t)=\prod_{j=0}^{n-1} f(d^j t) .\]
The problem is to find the behavior of the sequence of integrals
\[ I_n(f)=\int_0^1 f_n(t)\,dt \]
as $n\to\infty$. In particular, we want to find $c(f)$ defined by
\begin{equation}\label{L:c}
 c(f)=\limsup_{n\to\infty} |I_n(f)|^{1/n}.
 \end{equation}

The sequence $I_n$ is related to the bounded linear operator
\begin{equation}\label{T}
 (Tx)(t)=f(t)x(d\cdot t)
\end{equation}
which maps $L^2(\T)$ to itself. Note that
\[ f_n=T^n 1 \]
and
\begin{equation}\label{Im}
 I_n(f) =\langle T^n1,1\rangle
\end{equation}
with the inner product $\langle\cdot,\cdot\rangle$ in $L^2(\T)$.
In particular, 
\[ |I_n|\le \|T^n\|\]
and
\begin{equation}\label{L:eq1}
c(f)\le \lim_{n\to\infty} \|T^n\|^{1/n}=\rho(T).
\end{equation}

We show that $c(f)$ is equal to the spectral radius of a transfer operator under suitable assumptions on $f$. See also \cite{FanLau1998} for related results.

\begin{thm}\label{A:t1}
Suppose that $0\le f\in \H$ for some $0<\alpha\le 1$, and that the transfer operator $L$ defined by \eqref{L:L} is a Krein operator. Then $r=c(f)=R=\rho(L)$, where
$r,R$ are defined in \eqref{L:rR}. Moreover, we can replace $\limsup$ by $\lim$ in definition \eqref{L:c}.
\end{thm}
\begin{proof}
The adjoint $T^\ast$ of $T$ agrees with the operator $L$ when considered as an operator on $L^2(\T)$.
Let $h_n, r_n, R_n$ be defined by \eqref{L:hn}, \eqref{L:rnRn}.
It follows from \eqref{Im} that
\[ I_n=\langle 1, h_n\rangle=\int_0^1 h_n(t)\,dt .\]
Thus
\[r_n^{1/n}\le I_n^{1/n}\le R_n^{1/n}.\]
By Theorem \ref{L:t3}(e), the sequences $r_n^{1/n}$ and $R_n^{1/n}$ converge to the same limit $r=R$.
Therefore, the sequence $I_n^{1/n}$ converges and we obtain $r=c(f)=R$.
\end{proof}

Using Theorem \ref{A:t1} in connection with \eqref{L:est1} or Theorem \ref{L:t2} we can estimate $c(f)$. We will look at some examples in the next section.

We mention two special classes of functions $f$ for which $c(f)$ can be calculated explicitly.

1)
Suppose that $f$ is a step function such that $f(t)=f_i={\rm const}$ for $\frac{i-1}{d}\le t< \frac{i}{d}$, $i=1,\cdots,d$.
Then it is easy to show that
\[ I_n(f)= \left(\int_0^1 f(t)\,dt\right)^n.\]
Therefore,
\[ c(f)=\left|\int_0^1 f(t)\,dt\right| .\]
If $f$ is any nonnegative bounded measurable $1$-periodic function, we may introduce
two step function $g,h$ defined by
\begin{eqnarray*}
 g(t)&=&\inf\left \{ f(s): \frac{i-1}{d}\le s<\frac{i}{d}\right \}
 \quad\text{for $\frac{i-1}{d}\le t<\frac{i}{d}$},\\
h(t)&=&\sup\left \{ f(s): \frac{i-1}{d}\le s<\frac{i}{d}\right \}
\quad\text{for $\frac{i-1}{d}\le t<\frac{i}{d}$}.
 \end{eqnarray*}
Then we obtain the estimate
\begin{equation}\label{ineq}
I_n(g)\le I_n(f)\le I_n(h)
\end{equation}
and so
\begin{equation}\label{ineq1}
\int_0^1 g(t)\,dt\le c(f)\le \int_0^1 h(t)\,dt .
\end{equation}

2) Let $f$ be any bounded measurable $1$-periodic function with
Fourier expansion
\[ f(t)=\sum_{k\in\Z} a_k e^{2\pi i k t} .\]
We represent the operator $T$ by an infinite matrix in the orthonormal basis $\{e^{2\pi ikt}\}_{k\in\Z}$.
The matrix of $T$ is
\begin{equation}\label{T:matrix}
 \left( a_{k-d\ell} \right)_{k,\ell \in\Z} .
\end{equation}
In this notation $k$ is the row index and $\ell$ is the column index.
If we write
\[ f_n(t)=\sum_{k\in\Z} a_{k,n} e^{2\pi i kt},\]
then we obtain the coefficients $a_{k,n+1}$ from $a_{k,n}$ by application of $T$, so
\[ a_{k,n+1}=\sum_{\ell\in\Z} a_{k-d\ell}\, a_{\ell,n} .\]
Note that $I_n=a_{0,n}$.

In particular, suppose that $f(t)$ is a trigonometric polynomial of degree $N$, so
\[ f(t)=\sum_{k=-N}^N a_k e^{2\pi i k t} .\]
We set
\[ K:=\left\lfloor \frac{N-1}{d-1}\right\rfloor.\]
Consider the central $2K+1$ by $2K+1$ submatrix $B$ of $T$ consisting of rows $-K\le k\le K$ and columns $-K\le \ell\le K$.
Notice that all entries in the rows $-N\le k\le N$ outside the central submatrix vanish. Therefore, we obtain the recursion
\[ a_{k,n+1} =\sum_{\ell=-K}^K a_{k-d\ell}\,a_{\ell,n} \quad\text{if $-K\le k\le K$.}\]
Hence we can calculate $I_n=a_{0,n}$ by computing the powers of the matrix $B$.
It is clear that
\begin{equation}\label{eq2}
 c(f)\le \rho(B)
\end{equation}
but it is not immediately clear whether we have equality in \eqref{eq2}.
It depends on how the constant function $1$ is represented in a Jordan basis of $B$ (whether the basis vectors associated with largest eigenvalue of $B$ contribute to the expansion of $1$.)

The situation is clear if the matrix $B$ is nonnegative and primitive ($B^p$ is a positive matrix for some $p\in\N$.)
Then the spectral radius of $B$ is a simple positive eigenvalue and we can use Theorem 8.5.1 in \cite{HJ} to show that there is equality in \eqref{eq2}.
Suppose that $a_k>0$ for all $k=-N,-N+1,\cdots,N$. Then all entries in the main diagonal, the subdiagonal and superdiagonal of $B$ are positive. Therefore, $B$ is primitive.

If we have symmetry $a_{-k}=a_k$ then we can replace the matrix $B$ by a $K+1$ by $K+1$ matrix $C$ whose entries are
\[ c_{i,0}=a_i, \quad c_{i,j}=a_{i-dj}+a_{i+dj} \quad\text{if $0\le i\le K$, $1\le j\le K$.}\]
See the next section for examples.

\section{The special cases $f(t)=|\cos(\pi t)|^\q$ and $f(t)=|\sin(\pi t)|^\q$}\label{sec:special-weights}

In this section we consider the functions
\[f(t)=|\cos(\pi t)|^\q,\quad \tilde f(t)=|\sin(\pi t)|^\q\quad \text{where $\q>0$}.
\]
We set
\[ c(\q):=c(f),\quad \tilde c(\q):=c(\tilde f).\]
Obviously, $0\le c(\q), \tilde c(\q)\le 1$. Note that $f, \tilde f\in \H$ with $\alpha=\min(\q,1)$. By Lemma \ref{L:Krein2}, $L$ is a Krein operator except when $f(t)=|\cos(\pi t)|^\q$ and $d=2$.
This is an exceptional case that will be considered in the next section.
Except for this case we can apply Theorems \ref{L:t3} and \ref{A:t1}.

If $d$ is odd, then $I_n(f)=I_n(\tilde f)$, and consequently $c(f)=c(\tilde f)$. In fact, in this case the transfer operators weighted by $f$ and $\tilde f$ are conjugate to each other.

\begin{thm}\label{S:t1}
The functions $c(\q)$ and $\tilde c(\q)$ are convex and nonincreasing in $\q>0$. Moreover, 
\[ \lim_{\q\to\infty} c(\q)=\frac1 d, \quad
\lim_{\q\to\infty} \tilde c(\q)=\begin{cases} \frac1d & \text{if $d$ is odd,}\\ 0 & \text{if $d$ is even,}\end{cases}
\]
and
\[ \lim_{\q\to 0^+} c(\q)= \lim_{\q\to 0^+} \tilde c(\q)=1.\]
\end{thm}
\begin{figure}[h]
\includegraphics[scale=0.4]{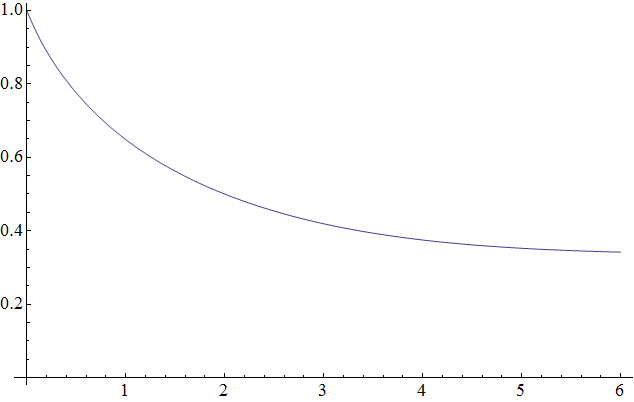}
\caption{A graph of $c(\q)$ for $0<\q<6$, when $d=3$.}
\vspace{2em}
\includegraphics[scale=0.4]{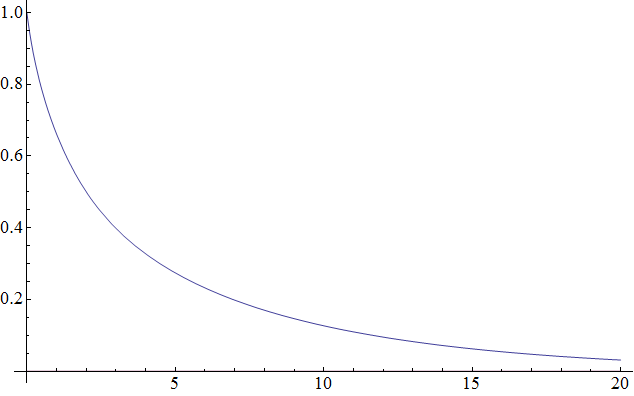}
\caption{A graph of $\tilde c(\q)$ for $0<\q<20$, when $d=2$.}
\end{figure}
\begin{proof}
Monotonicity of $c(\q)$ and $\tilde c(\q)$ are clear. Convexity follows from H\"older's inequality and Young's inequality. By looking at the function $h_1$ associated with $f$, we see that $R_1=\max_{t\in[0,1]} h_1(t)$ converges to $d^{-1}$ as $\q\to\infty$. So
\[\limsup_{q\to\infty} c(\q)\le d^{-1}.\]
If we use an estimate of the form $f(t)\ge 1-\frac{t}{a}$ for $0\le t\le a$ for some $a\in(0,1)$ depending on $\q$, then we can estimate
\[ f_n(t)\ge \prod_{j=0}^{n-1} \left(1-\frac{d^j t}{a}\right)\ge \exp\left(-(\ln 4)\, \frac{d^n-1}{d-1} \frac{t}{a}\right)\quad\text{for $0\le t\le \frac{a}{2 d^{n-1}}$} .\]
It then follows from \eqref{L:c} that $c(\q) \ge d^{-1}$ for all $q>0$. Hence $$\lim_{\q\to\infty} c(\q)=d^{-1}.$$

When treating $\tilde c(\q)$, we may assume that $d$ is even. The proof is similar to the preceding one. To determine the limit of $\tilde c(\q)$ as $\q\to\infty$, we use $c(\tilde f)\le R_2^{1/2}$.

To show the limits as $\q\to 0^+$, it suffices to show that
\[c(\q)^{1/\q}\ge 1/2,\quad \tilde c(\q)^{1/\q}\ge 1/2\]
for all $q>0$. To this end,
let $g(t)=|\cos(\pi t)|$ (or $|\sin(\pi t)|$). By Jensen's inequality, we have
$$I_n^{1/\q}=\left (\int_0^1 |g_n(t)|^\q dt\right )^{1/\q} 
\ge \exp\left (\int_0^1 \ln |g_n(t)|dt\right )$$
for all $\q>0$. 
On the other hand, since
$$\int_0^1 \ln |\cos(\pi t)|dt=\int_0^1 \ln |\sin(\pi t)|dt=\ln (1/2),$$
we have
$$\int_0^1 \ln |g_n(t)|dt=\sum_{j=0}^{n-1} \int_0^1 \ln |g(d^j t)|dt=n\ln (1/2).$$
Thus, for any $\q>0$,
$$c(\q)^{1/\q}\ (\text{or }\tilde c(\q)^{1/\q})=\lim_{n\to\infty} I_n^{1/(n\q)}\ge 1/2.$$
This completes the proof.
\end{proof}

\subsection{The case $f(t)=\cos^{2N}(\pi t)$, $d=3$}
For $N\in\N$, consider the trigonometric polynomial
\[ f(t)=\cos^{2N}(\pi t). \]
The degree of $f$ is $N$ and, for $-N\le k\le N$,
\[ a_k=2^{-2N} \binom{2N}{k+N}>0 .\]
We use the method 2) from Section \ref{sec:quasicompact}.
If $N=1,2$, then $K=0$ and so $c(2)=\frac12$ and $c(4)=\frac38$ (note that this recovers a result of Strichartz \cite{Strichartz1990}).
If $N=3$, then $K=1$ and
\[ C=\frac1{64} \left[\begin{array}{rr} 20 & 2\\ 15 & 6 \end{array}\right]. \]
So
\[ c(6)=\rho(C)=\frac1{64}(13+\sqrt{79}) =0.342003\cdots.\]
If $N=4$, then $K=1$ and
\[ C=\frac1{256} \left[\begin{array}{rr} 70 & 16\\ 56 & 29 \end{array}\right]. \]
It follows that
\[ c(8)=\rho(C)= \frac1{512}(99+9\sqrt{65})=0.335078\cdots.\]
When $N\ge 5$, the formulas for $c(2N)$ become more complicated, but $c(2N)$ is easy to compute numerically.
For example, we obtain
\[ c(10)=0.333691\cdots.\]

The same method can be used to determine $c(2N)$ for other values of $d$.

\subsection{The case $f(t)=|\sin(\pi t)|$, $d=3$}\label{subsec:d=3}
Even if $f$ is not a trigonometric polynomial, we can still use matrix methods to estimate $c(f)$.
As an example, consider
\begin{equation}\label{S:fs}
 f(t)=|\sin(\pi t)|.
\end{equation}
The Fourier coefficients of $f$ are
\[ a_k=-\frac{2}{\pi} \frac{1}{(2k-1)(2k+1)},\ k\in\Z .\]
Note that $a_0>0$ but all other $a_k$ are negative.
Let $N\in\N$. We estimate
\[ f(t)\le  \sum_{k=-N}^N a_k e^{2\pi i k t} +\frac2\pi \sum_{k=N+1}^\infty \frac{2}{(2k-1)(2k+1)} \]
so we have
\[ f(t)\le h(t) ,\]
where
\[ h(t)=\frac2\pi \frac1{2N+1} +  \sum_{k=-N}^N a_k e^{2\pi i k t} .\]
Using
\[ c(f)\le c(h)\le \rho(C) \]
we get upper bounds for $c(f)=c(1)$ (when $d=3$):

\begin{center}
\begin{tabular}{r|c}
$N$ & $\rho(C)$ \\
\hline
1 & 0.848826\\
2 & 0.763943\\
3 & 0.737463\\
4 & 0.717381\\
5 & 0.704696\\
10& 0.678384\\
20& 0.663593\\
30& 0.658613\\
50 & 0.654552\\
100 & 0.651436
\end{tabular}
\end{center}
As far as we know the exact value of $c(1)$ is not known. We conjecture that
$c(1)=0.648314\cdots$.

We also obtain
\[ g(t)\le f(t) \]
where
\[ g(t)= -\frac2\pi \frac1{2N+1} + \sum_{k=-N}^N a_k e^{2\pi i k t} .\]
Since $a_0>0$ and all other $a_k<0$, we can easily show that $g(t)\ge 0$.
Therefore, we have
\[ c(g)\le c(f) .\]
However, the trigonometric polynomial $g$ does not have positive coefficients, so we do not know whether $c(g)=\rho(C)$. Therefore, we do not obtain lower bounds by this method.
For $N=100$, one would get $\rho(C)=0.645194\cdots $.

Somewhat surprisingly, the functions $h_n$ associated with \eqref{S:fs}
can be represented in a fairly explicit way.
If
\[ u(t)=\cos\left(\pi as\right),\quad s=t-\tfrac12, \]
then
\[ (Lu)(t)= \tfrac16 \big(1+2\cos\tfrac\pi3(1+a)\big)\cos\tfrac\pi3(1+a)s+\tfrac16\big(1+2\cos\tfrac\pi3(1-a)\big) \cos\tfrac\pi3(1-a) s .\]
Iterating this formula, we see that $h_n$ is a sum of $2^{n-1}$ many terms of the form
\[ A\cos(\pi a s) \]
where $A>0$ and $a=\frac13\pm\frac19\pm\frac1{27}\pm\cdots\pm \frac{1}{3^n}\in (0,\tfrac12)$.
It follows that
\[ r_n=h_n(0),\quad R_n=h_n(\tfrac12).\]
By \eqref{L:rR},
\[ \left(h_n(0)\right)^{1/n}\le c(f)\le \left(h_n(\tfrac12)\right)^{1/n} .\]
For example, if $n=1$ we get the bounds
\[ \frac13\sqrt3\le c(f)\le \frac23.\]

Since $0<a<\frac12$, we have
\[ R_n\le \sqrt{2} r_n. \]
Therefore,
\[ (R_n)^{1/n} \le 2^{1/(2n)} (r_n)^{1/n} .\]
In agreement with Theorem \ref{A:t1}, we get
\[ r=R=c(f) .\]
We also find that
\[ (R_n)^{1/n}-(r_n)^{1/n}\le (R_n)^{1/n}(1-2^{-1/(2n)})\le \tfrac23 (1-2^{-1/(2n)}) .\]
For example, if $n=10$, then  $c(f)$ is enclosed in the interval $[(r_n)^{1/n},(R_n)^{1/n}]$ of length at most $\frac23(1-2^{-1/20})=0.0277\cdots$ (the actual length is
$0.008390\cdots$.)

By numerical computation, we get the following bounds:
\begin{center}
\begin{tabular}{r|r|r}
$n$ & $(h_n(0))^{1/n}$ & $(h_n(\tfrac12))^{1/n}$ \\
\hline
1 & 0.577350 & 0.666666  \\
2 & 0.615672 & 0.656538  \\
3 & 0.626102 & 0.653844  \\
4 & 0.631603 & 0.652453  \\
5 & 0.634908 & 0.651623   \\
10& 0.641576 & 0.649967  \\
15& 0.643815 & 0.649415
\end{tabular}
\end{center}

\subsection{Properties of $h_n(t)$ when $f(t)=|\sin(\pi t)|^\q$}
If we use $r_n^{1/n}$ and $R_n^{1/n}$ to bound $c(f)$, we are faced with the problem to compute the maximum and minimum values of the function $h_n$.
Therefore, it is of interest to discuss the behavior of the function $h_n$.
Consider
\[ f(t)=|\sin(\pi t)|^\q, \quad d\ge 2.\]
Then we have
\[ h_1(t)=\frac1d \sum_{i=0}^{d-1} \left|\sin\Big(\pi \frac{t+i}{d}\Big)\right|^\q,\quad 0\le t\le 1. \]
Note that $h_1(t)=h_1(1-t)$.
For this function we have the following result.

\begin{lemma}\label{S:l1}
(a) If $\q=2,4,6,\cdots,2(d-1)$, then 
$$h_1(t)\equiv \frac{1}{2^\q}\binom{\q}{\q/2}.$$\\
(b) Define the intervals 
\[ Q_k=\begin{cases} \big(2(k-1),2k\big) & \text{if $k=1,\cdots,d-1$,}\\ 
\big(2(d-1),\infty\big) & \text{if $k=d$.}\end{cases}
\]
Then for $\q\in Q_k$ we have
$$(-1)^{k-1}h_1'(t)>0,\quad 0<t<1/2.$$
\end{lemma}

\begin{proof}
(a) follows from the matrix representation of $T^\ast$, the adjoint of the matrix \eqref{T:matrix}.\\
(b) We differentiate $h_1(t)$ to get
\begin{equation}\label{diff}
 \frac{d^2}{\pi q} h_1'(t)=\sum_{j=1}^d s_j^\q c_j,\quad 0<t<1,
 \end{equation}
where
\[ s_j=\sin\left(\pi\frac{t+j-1}{d}\right),\quad c_j=\cot\left(\pi\frac{t+j-1}{d}\right) .\]
By (a), the left-hand side of \eqref{diff} is zero for $\q=2,4,\cdots,2(d-1)$.
Therefore, we obtain a linear system $V c= b$, where $V$ is a $d\times d$ generalized Vandermonde matrix with entries $s_j^{\q_i}$ in the $i$th row and $j$th column, where $\q_1=\q$ and $\q_i=2(i-1)$ for $i=2,3,\cdots,d$. The column vector $c$ has components $c_i$, and the column vector $b$
has first component $\frac{d^2}{\pi q} h_1'(t)$ and all other components equal to $0$.
Suppose that $0<t<\frac12$ and $\q\in Q_k$ for some $k=1,2,\cdots,d$.
Then both $s_1,s_2,\cdots,s_d$ and $\q_1,\q_2,\cdots,\q_d$ are mutually distinct. It is known (cf.\cite[page~76]{GK}) that
this implies that $\det V\ne 0$, and  $\det V>0$ if $s_1<s_2<\cdots<s_d$ and $\q_1<\q_2<\cdots<\q_d$.
If we solve the linear system $Vc=b$ for $c$ by Cramer's rule we find that $c_1$ has the same sign as $(-1)^{k-1} h_1'(t)$.
Since $c_1>0$, we have $(-1)^{k-1}h_1'(t)>0$ when $0<t<\frac12$ and $\q\in Q_k$.
\end{proof}

Using the bounds $r_1\le c(f)\le R_1$ together with Lemma \ref{S:l1}, one can deduce for $d=3$
$$\begin{cases}
\frac23\left(\frac{\sqrt3}{2}\right)^\q \le  c(\q) \le  \frac13+\frac23\left(\frac12\right)^\q&\quad \text{if $0<\q\le 2$ or $\q\ge 4$},\\
\frac13+\frac23\left(\frac12\right)^\q\le  c(\q) \le \frac23\left(\frac{\sqrt3}{2}\right)^\q&\quad \text{if $2\le \q\le 4$}.
\end{cases}$$
Note also that Lemma \ref{S:l1}(a) implies
$$c(q)=\tilde c(q)= \frac{1}{2^\q}\binom{\q}{\q/2}\quad \text{when } q=2,4,6,\cdots,2(d-1).$$

We would like to extend Lemma \ref{S:l1} to $h_2, h_3,\cdots$.
Based on computer experiments we conjecture the following.

\begin{conj}\label{S:conj1}
Lemma \ref{S:l1}(b) is true for every $h_n$, $n\in\N$.
\end{conj}

We obtain sharper lower and upper bounds for $c(f)$ when we
choose $w(t)=h_n(t)$ in \eqref{L:est2}. More precisely, we get
\begin{equation}\label{S:q3}
\min_{t\in[0,1]} \frac{h_{n+1}(t)}{h_n(t)}\le c(f)\le \max_{t\in[0,1]} \frac{h_{n+1}(t)}{h_n(t)} .
\end{equation}
Here we are faced with the problem
to determine the extrema of the quotients $\frac{h_{n+1}(t)}{h_n(t)}$.
Computer calculations suggest the following.

\begin{conj}\label{S:conj2}
Let $d=3$. If $0<\q<2$ and $n$ is odd, then $\frac{h_{n+1}(t)}{h_n(t)}$ attains its maximum at $t=0$ and its minimum at $t=\frac12$.
If $0<\q<2$ and $n$ is even, then $\frac{h_{n+1}(t)}{h_n(t)}$ attains its maximum at $t=\frac12$ and its minimum at $t=0$.
If $2<\q<4$, then $\frac{h_{n+1}(t)}{h_n(t)}$ attains its maximum at $t=0$  and its minimum at $t=\frac12$.
If $\q>4$, then $\frac{h_{n+1}(t)}{h_n(t)}$ attains its maximum at $t=\frac12$  and its minimum at $t=0$.
\end{conj}

If we believe these conjectures then $c(f)$ would lie between $\frac{r_{n+1}}{r_n}$ and $\frac{R_{n+1}}{R_n}$ for every $n$.
In the case $d=3$ we get the following estimates for $c(1)$:

\begin{center}
\begin{tabular}{r|r|r}
$n$ & lower bound & upper bound \\
\hline
1 & 0.577350 & 0.666666  \\
2 & 0.646564 & 0.656538  \\
3 & 0.648297 & 0.648396
\end{tabular}
\end{center}

We see that these bounds are much better than those from Section \ref{subsec:d=3}. Unfortunately, we used conjecture \ref{S:conj2} but for small $n$ it can be proved by direct computation.

\section{The case $f(t)=|\cos(\pi t)|^\q$, $d=2$}\label{sec:binary}

In the exceptional case
\[ f(t)=|\cos(\pi t)|^\q, \quad d=2\]
we can obtain more explicit computations. See also \cite{FanLau1998} for related results.

\subsection{The integrals $I_n$}
Using the identity
\[
\prod_{j=0}^{n-1} \cos(2^j t)=\frac{\sin(2^n t)}{2^n \sin(t)} ,
\]
we can write
\begin{equation}\label{1:fm}
f_n(t)=\frac{1}{2^{\q n}} \frac{|\sin(\pi 2^n t)|^\q}{|\sin(\pi t)|^\q}.
\end{equation}
Therefore,
\begin{equation}\label{1:fm-integral}
I_n=\frac{1}{2^{\q n}} \int_0^1 \frac{|\sin(\pi 2^n t)|^\q}{|\sin(\pi t)|^\q}\, dt=\frac{1}{2^{\q n-1}} \int_0^{1/2} \frac{|\sin(\pi 2^n t)|^\q}{|\sin(\pi t)|^\q}\, dt.
\end{equation}

\begin{thm}\label{1:t}
For $d=2$, we have
\[ c(\q)=\lim_{n\to\infty} I_n^{1/n}=
\begin{cases} 
2^{-\q} & \text{if $0<\q\le 1$,}\\
\frac12 & \text{if $\q>1$.}
\end{cases}
\]
\end{thm}
\begin{proof}
Substituting $u=2^n t$ in \eqref{1:fm-integral}, we get
\[ I_n=\frac{1}{2^{n-1}} \int_0^{2^{n-1}} \frac{|\sin(\pi u)|^\q}{2^{\q n}|\sin(\pi 2^{-n} u)|^\q}\,du .\]
Using
\[ \frac{2}{\pi} t\le \sin t\le t,\quad 0\le t\le \frac{\pi}{2},\]
we find
\begin{equation}\label{1:ineq}
\frac{\pi^{-\q}}{2^{n-1}} \int_0^{2^{n-1}} \frac{|\sin(\pi  u)|^\q}{u^\q}\, du\le I_n \le \frac{2^{-\q}}{2^{n-1}} \int_0^{2^{n-1}} \frac{|\sin(\pi  u)|^\q}{u^\q}\, du  .
\end{equation}
If $q>1$, the integral
\begin{equation}\label{1:int}
\int_0^\infty \frac{|\sin(\pi  u)|^\q}{u^\q}\,du
\end{equation}
converges. Therefore, the statement of the theorem follows for $\q>1$.
If $\q=1$, the integral \eqref{1:int} diverges and the integrals in \eqref{1:ineq} behave like $\ln(2^n)$. Since $n^{1/n}$ converges to $1$ as $n\to\infty$,
we obtain the statement of the theorem when $\q=1$. If $0<\q<1$, the integrals in \eqref{1:ineq}  behave like $2^{n(1-\q)}$ which implies the statement of the theorem
for $0<\q<1$.
\end{proof}

\subsection{Spectral radius} 
By \eqref{L:Ln}, we have
\[ h_n(t)=\frac{1}{2^{n}}\sum_{k=0}^{2^n-1} f_n\left (\frac{t+k}{2^{n}}\right ) .\]
Combining with \eqref{1:fm}, we get
\begin{equation}\label{2:hm}
h_n(t)=\frac{|\sin(\pi t)|^\q}{2^n 2^{\q n}} \sum_{k=0}^{2^n-1} \frac{1}{|\sin(\pi2^{-n}(t+k))|^\q}.
\end{equation}
By estimating the sum in \eqref{2:hm} we obtain the following.

\begin{thm}\label{2:t}
For $f(t)=|\cos(\pi t)|^\q$, $d=2$, we have
\[ r=\frac12 \quad  \text{for all $\q>0$}\]
and
\[ R=\begin{cases}
2^{-\q} & \text{if $0<\q\le 1$,}\\
\frac12 & \text{if $\q>1$.}
\end{cases}
\]
\end{thm}
\begin{proof}
Using only the term with $k=0$ in \eqref{2:hm}, we obtain, for $0\le t\le \frac12$,
\[ h_n(t)\ge \frac{1}{2^n}\frac{1}{2^{\q n}}\frac{|\sin(\pi t)|^\q}{|\sin(\pi 2^{-n}t)|^\q}\ge \frac{1}{2^n}\frac{2^\q}{\pi ^\q}.\]
This inequality together with $h_n(0)=2^{-n}$ proves $r=\frac12$.
The proof of the formula for $R$ is elaborated in Section \ref{subsec:convergence}.
\end{proof}

\subsection{Eigenfunctions}\label{subsec:eigenfunctions}
Let $\alpha=\min\{1,\q\}$.
By Theorem \ref{L:t1}, the spectral radii of $L$ and $L_\alpha$ agree, and $L_\alpha$ is quasicompact.
Since $L$ is also a positive operator, $\lambda=R$ must be an eigenvalue, so there must exist a corresponding eigenfunction.
But $L$ is not a Krein operator (cf. \cite{A}), so we do not know whether the eigenfunction is unique (up to a constant factor) 
or whether it is positive on $\T$.

We want to find nontrivial solutions $u\in C(\T)$ to the equation  $Lu=\lambda u$, particularly for $\lambda=R$.
Interestingly, we can find these eigenfunctions fairly explicitly.
In fact, if we substitute
\[ u(t)=|\sin(\pi t)|^\q g(t) \]
in $Lu=\lambda u$, we find
\begin{equation}\label{3:ber}
\tfrac12 g\left(\tfrac12 t\right)+\tfrac12 g\left(\tfrac12 (t+1)\right) =\mu g(t)\quad \text{where }\mu=2^\q \lambda .
\end{equation}
Note that $g(t)$ will usually be continuous only on the open interval $(0,1)$.
Much is known about equation \eqref{3:ber} (cf. \cite{V}).
Clearly, $g(t)=1$ is a solution to \eqref{3:ber} with $\mu=1$. Therefore, $u(t)=|\sin(\pi t)|^\q$ is an eigenfunction of $L$ corresponding
to the eigenvalue $\lambda=2^{-\q}$. If $0<\q\le 1$, then this is an eigenfunction corresponding to the spectral radius eigenvalue $R$.
Furthermore, $g(t)=B_n(t)$ with $B_n(t)$ denoting a Bernoulli polynomial, is also a solution to \eqref{3:ber} corresponding to $\mu=2^{-n}$.
This gives us many more eigenfunctions of $L$, but they do not give us eigenfunctions corresponding to the eigenvalue $\lambda=R$ if $\q>1$.

Using an idea from \cite{V}, we find eigenfunctions corresponding to $\lambda=R$ when $\q>1$.
For $s>1$, consider the Hurwitz zeta function
\begin{equation}\label{3:g}
\zeta(s,t)=\sum_{k=0}^\infty \frac{1}{(t+k)^s},\quad 0<t<1.
\end{equation}
It is easy to check that $g(t)=\zeta(s,t)$ is a solution to \eqref{3:ber} with $\mu =2^{s-1}$.
If we let 
\begin{equation}\label{3:Gst}
G(t)=G(s,t):=g(t)+g(1-t),\quad 0<t<1.
\end{equation}
Then $G(t)$ is also a solution to \eqref{3:ber} and has symmetry $G(t)=G(1-t)$. If $s\le \q$, 
then
\[u(t)=|\sin(\pi t)|^\q G(t)\]
is a continuous eigenfunction of $L$ corresponding to the eigenvalue $\lambda=2^{s-\q-1}$.
In particular, choosing $s=\q$, we obtain an eigenfunction corresponding to the eigenvalue $R=\frac12$.

Suppose $\q\ge 2$ is an even integer. Let $s=2,3,\cdots,\q$.
Consider 
\[ \widetilde G(t)=\zeta(s,t)+(-1)^s\zeta(s,1-t)=\sum_{k=-\infty}^\infty \frac{1}{(t+k)^s}. \]
Since
\[ \pi\cot(\pi t)=\lim_{N\to\infty} \sum_{k=-N}^N \frac1{t+k}, \]
we obtain
\begin{equation}\label{3:G}
\widetilde G(t)=\frac{(-1)^{s-1}\pi}{(s-1)!}\left(\frac{d}{dt}\right)^{s-1} \cot(\pi t) ,
\end{equation}
and correspondingly the eigenfunctions 
\[ u(t)=|\sin(\pi t)|^\q \widetilde G(t) .\]
Obviously, these eigenfunctions are trigonometric polynomials. For example, if $\q=s=2$, we obtain the eigenfunction
\[ u(t)=|\sin(\pi t)|^2 \widetilde G(t)=|\sin(\pi t)|^2 G(t)= \pi^2 ,\]
and, if $\q=s=4$, then
\begin{equation}\label{3:q4}
u(t)=\frac{\pi^4}{3}\big(\cos(2\pi t)+2\big) .
\end{equation}

We can also find these eigenfunctions in a different way.
The space of trigonometric polynomials $\sum_{k=-K}^K c_k e^{2\pi i kt}$ with $K=(\q-2)/2$ is an invariant subspace of $L$.
The matrix representation of the restriction of $L$ to this invariant subspace with respect to the basis $\{e^{2\pi i kt}\}$ is
\[ B=(a_{2\ell-k})_{-K\le k, \ell\le K},\]
where the $a_k$ denote the Fourier coefficients in
\[ |\cos(\pi t)|^\q= \sum_{k=-K}^K a_k e^{2\pi ikt}.\]
For example, if $\q=4$ then
\[ B=\left[\begin{array}{rrr} 
\frac14 & \frac{1}{16} & 0 \\[0.3em]
\frac14 & \frac38 & \frac14 \\[0.3em]
0 & \frac1{16} & \frac14 \\
\end{array}\right]
\]
This matrix is nonnegative and primitive. Its largest eigenvalue is $\frac12$, which follows from the fact that the column sums are all equal to $\frac12$.
An eigenvector corresponding to the eigenvalue $\lambda=\frac12$ is $[1,4,1]^T$.
Therefore,
\[ u(t)=\cos(2\pi t)+2\]
is an eigenfunction of $L$ for $\q=4$ corresponding to the eigenvalue $\lambda=\frac12$.
Apart from a constant factor, this is the same eigenfunction we found in \eqref{3:q4}.
This is true in general. The eigenfunctions of the form $u(t)=|\sin(\pi t)|^\q \widetilde G(t)$ with $\widetilde G$ from \eqref{3:G} with $s=2,3,\cdots,q$ match the eigenfunctions obtained from the matrix $B$. In particular, we see that the matrix $B$ has eigenvalues $\frac12, \frac14,\cdots,2^{-\q+1}$.

In Section \ref{subsec:h_n-properties}, we study in more detail the trigonometric polynomials obtained from \eqref{3:G}.

\subsection{Convergence of $h_n(t)$}\label{subsec:convergence}
The following proposition shows that, after appropriately normalized, the function $h_n$ converges to the eigenfunctions we found in Section \ref{subsec:eigenfunctions} corresponding to $\lambda=R$.

\begin{prop}\label{prop:convergence} The following limits hold in $C^\q(\T)$.\\
(a) If $0<\q<1$, then
$$\lim_{n\rightarrow\infty} 2^{\q n} h_n(t)=\frac{\Gamma(\frac{1}{2}-\frac{\q}{2})}{\sqrt{\pi}\Gamma(1-\frac{\q}{2})}|\sin(\pi t)|^\q.$$
(b) If $\q=1$, then
$$\lim_{n\rightarrow\infty} \frac{2^{n}}{n} h_n(t)=\frac{2\ln 2}{\pi}|\sin(\pi t)|.$$
(c) If $1<\q<\infty$, then
$$\lim_{n\rightarrow\infty} 2^{n} h_n(t)=\pi^{-\q}|\sin(\pi t)|^\q G(\q,t).$$
where $G(\q,t)$ is given by \eqref{3:Gst}.
\end{prop}

\begin{proof}
We show here the pointwise convergence of $h_n$. The norm convergence can be shown by slight refinements of the argument.\\
(a) By \eqref{2:hm}, we have
$$h_n(t)=\frac{|\sin(\pi t)|^\q}{2^n 2^{\q n}} \sum_{k=0}^{2^n-1} \frac{1}{|\sin(\pi2^{-n}(t+k))|^\q}.$$
Since 
$$\int_0^1 \frac{1}{|\sin(\pi x)|^\q}dx=\frac{\Gamma(\frac{1}{2}-\frac{\q}{2})}{\sqrt{\pi}\Gamma(1-\frac{\q}{2})},$$
to prove the statement it suffices to show that
$$\lim_{n\rightarrow\infty}\frac{1}{2^{n}}\sum_{k=0}^{2^n-1} \frac{1}{|\sin(\pi 2^{-n}(t+k))|^\q}=\int_0^1 \frac{1}{|\sin(\pi x)|^\q}dx.$$
However, this follows easily by treating the left-hand side as a Riemann sum, using the monotonicity of the integrand and the assumption that $\q<1$.

\noindent(c) Similar as in the proof of (a), we only need to show that
$$\lim_{n\rightarrow\infty}\frac{1}{2^{\q n}}\sum_{k=0}^{2^n-1} \frac{1}{|\sin(\pi 2^{-n}(t+k))|^\q}=\pi^{-\q} G(\q,t).$$
By symmetry, this reduces to showing
$$\lim_{n\rightarrow\infty}\frac{1}{2^{\q n}}\sum_{k=0}^{2^{n-1}-1} \frac{1}{|\sin(\pi 2^{-n}(t+k))|^\q}=\pi^{-\q}\sum_{k=0}^{\infty} \frac{1}{(t+k)^\q}.$$
However, this follows easily from the basic limit
$$\lim_{x\rightarrow 0}\frac{\sin x}{x}=1$$
and the fact that the series on the right-hand side is convergent.

\noindent(b) By symmetry, it suffices to show that
$$\lim_{n\rightarrow\infty}\frac{1}{n 2^{n}}\sum_{k=0}^{2^{n-1}-1} \frac{1}{\sin(\pi 2^{-n}(t+k))}=\frac{\ln 2}{\pi}.$$
To this end, for any given $\varepsilon>0$ we fix $\delta>0$ such that
$$\left |\frac{x}{\sin x}-1\right |<\varepsilon,\quad 0<x<\delta.$$
We can then write
\begin{align*}
&\sum_{k=0}^{2^{n-1}-1} \frac{1}{\sin(\pi 2^{-n}(t+k))}\\
=&\sum_{k:\,\pi 2^{-n}(t+k)<\delta}\frac{\pi 2^{-n}(t+k)}{\sin(\pi 2^{-n}(t+k))}\frac{1}{\pi 2^{-n}(t+k)}\\
&+\sum_{k:\,\pi 2^{-n}(t+k)\ge\delta}\frac{1}{\sin(\pi 2^{-n}(t+k))}\\
=&I+II.
\end{align*}
By our choice of $\delta$,
\begin{align*}
I
&=(1+O(\varepsilon))\sum_{k:\,\pi 2^{-n}(t+k)<\delta} \frac{1}{\pi 2^{-n}(t+k)}\\
&=\frac{2^n}{\pi}(1+O(\varepsilon))\sum_{k:\,\pi 2^{-n}(t+k)<\delta} \frac{1}{t+k}\\
&=\frac{2^n}{\pi}(1+O(\varepsilon))(1+o(1))\ln (2^n)\\
&=n 2^n(1+O(\varepsilon))\frac{\ln 2}{\pi},\quad \text{as } n\to\infty.
\end{align*}
On the other hand, using
$$\sin x\ge \frac{2}{\pi}x,\quad 0\le x\le\pi/2,$$
we have
\begin{align*}
II
&=O(2^n) \sum_{k:\,\pi 2^{-n}(t+k)\ge\delta}\frac{1}{t+k}\\
&=O(2^n) \ln\Big(\frac{1}{\delta}\Big)\\
&=o(n 2^n),\quad \text{as } n\to\infty.
\end{align*}
Combining these, we get
$$\frac{1}{n 2^n}(I+II)=(1+O(\varepsilon))\frac{\ln 2}{\pi},\quad \text{as } n\to\infty.$$
Since $\varepsilon$ is arbitrary, this completes the proof.
\end{proof}

\subsection{Properties of $h_n(t)$}\label{subsec:h_n-properties} It turns out that the functions $h_n, n\in\N$ share some common geometric properties. We were able to prove some of them.

\begin{prop}\label{prop:convexity}
(a) If $0<\q\le 1$, then $h_n''(t)<0$ for all $n\in\N$ and $t\in (0,1)$.\\
(b) If $1<\q<2$, then $h_n''(1/2)<0$ for $n=1,\cdots,N(\q)$, where $N(\q)$ satisfies $\lim_{\q\rightarrow 1^+}N(\q)=\infty$.\\
(c) If $2<\q<\infty$, then $h_n''(1/2)>0$ for $n=1,\cdots,N(\q)$, where $N(\q)$ satisfies $\lim_{\q\rightarrow\infty}N(\q)=\infty$.
\end{prop}
\begin{proof}
\noindent(a) In the case $n=0$, $h_0(t)\equiv 1$, so the statement obviously holds with strict inequality replaced by equality. Assume that $h_{n-1}''(t)\le 0$ for all $t\in (0,1)$. We now show that $h_n''(t)<0$ for all $t\in (0,1)$.

By definition, we have
$$h_n(t)=\frac{1}{2}\Big|\cos\Big(\pi \frac{t}{2}\Big)\Big|^\q h_{n-1}\Big(\frac{t}{2}\Big)+\frac{1}{2}\Big|\cos\Big(\pi \frac{t+1}{2}\Big)\Big|^\q h_{n-1}\Big(\frac{t+1}{2}\Big).$$
Since the second term equals the first term after the change of variable $t\rightarrow 1-t$, it suffices to show $(fg)''(t)<0$ for all $t\in (0,1/2)$, where
$$f(t)=|\cos(\pi t)|^\q,\quad g(t)=h_{n-1}(t).$$
However, by the product rule,
$$(fg)''=f''g + 2 f'g' + fg''.$$
Since $\q\le 1$, we have $f''(t)<0$ for all $t\in (0,1/2)$. Also, by symmetry we have $g'(1/2)=0$, and so the induction hypothesis implies  $g'(t)\ge 0$ for all  $t\in (0,1/2)$. Combining these we get $f''g<0$, $f'g'\le 0$, and $fg''\le 0$, which gives $(fg)''(t)<0$ for all $t\in (0,1/2)$. This completes the proof by induction.

\noindent(b) The proof is similar to that of (a). Using the same notation, we observe that
$$f''(t)=\pi^2 \q |\cos(\pi t)|^\q \Big(\q|\sin(\pi t)|^2-1\Big),\quad 0<t<1/2$$
now changes sign at
$$t=t_\q:=\pi^{-1} \arcsin(\sqrt{1/\q}).$$
Notice that $t_\q>1/4$ if $\q<2$, and $t_\q<1/4$ if $\q>2$; moreover,
$$\lim_{\q\rightarrow 1^+} t_\q = 1/2,\quad
\lim_{\q\rightarrow \infty} t_\q = 0.$$
In order to determine the sign of
\begin{equation}\label{eqn:2nd-derivative}
(fg)''=f''g + 2 f'g' + fg'',
\end{equation}
as before we want all the three terms to have the same sign.

In the case $n=1$, since $g\equiv 1$, we have $(fg)''(t)=f''(t)<0$ for all $t\in(0,2t_\q)$, where $2t_\q>1/2$ and $2t_\q\rightarrow 1$ as $\q\rightarrow 1^+$. This implies $h_1''(t)<0$ for all $t\in (1-2t_\q,2t_\q)$. By symmetry we have $h_1'(t)>0$ for all $t\in (1-2t_\q,1/2)$. Now proceeding by induction, we see that, using \eqref{eqn:2nd-derivative},
$$h_n''(t)<0,\quad t\in (2^{n-1}(1-2t_\q),1-2^{n-1}(1-2t_\q))$$
and
$$h_n'(t)>0,\quad t\in (2^{n-1}(1-2t_\q),1/2).$$
In particular, if $\q>1$ is sufficiently close to 1, we have $2^{n-1}(1-2t_\q)<1/2$ and thus $h_n''(1/2)<0$, as desired.

The proof for (c) is similar.
\end{proof}

When $\q$ is an even integer, we can have more information.

\begin{prop}
If $\q\ge 4$ is an even integer, then
$$h_\infty(t):=\pi^{-\q}|\sin(\pi t)|^\q G(\q,t)$$
satisfies $h_\infty''(1/2)>0$; moreover, $h_\infty'(t)<0$ for all $t\in(0,1/2)$.
\end{prop}

\begin{proof}
By \eqref{3:G}, we have
$$G(\q,t)=\frac{(-1)^{\q-1}\pi}{(\q-1)!}
\left (\frac{d}{dt}\right )^{\q-1}\cot (\pi t).$$
Lemma \ref{lem:polynomial} below shows that, after simplification,
$$h_\infty(t)=\frac{2}{(\q-1)!} P_{\q-1}(\cos(\pi t))$$
where $P_{\q-1}(x)$ is a polynomial consisting of the even powers $1, x^2, \cdots, x^{\q-2}$ and has positive coefficients. By direct computation, we then have
$$h_\infty'(t)=-\frac{2\pi}{(\q-1)!}P_{\q-1}'(\cos(\pi t)) \sin(\pi t)$$
and
$$h_\infty''(1/2)=\frac{2\pi^2}{(\q-1)!} P_{\q-1}''(0).$$
The desired conclusions now follow immediately from the properties of $P_{\q-1}(x)$ mentioned above.
\end{proof}

\begin{lemma}\label{lem:polynomial}
For all $n\in\N$, we have
$$\left (\frac{d}{dt}\right )^{n}\cot t=(-1)^n\frac{P_n(\cos t)}{(\sin t)^{n+1}}$$
where $P_n(x)$ is a polynomial of degree $n-1$ whose coefficients are nonnegative integers. Moreover, when $n$ is odd, $P_n(x)$ consists of the even powers $1,x^2,\cdots,x^{n-1}$; when $n$ is $even$, $P_n(x)$ consists of the odd powers $x,x^3,\cdots,x^{n-1}$.
\end{lemma}

\begin{proof}
It is easy to see that $P_0(x)=x$ and $P_1(x)=1$. Moreover, by direct computation we have
$$P_{n+1}(x)=(n+1)xP_{n}(x)+(1-x^2)P_{n}'(x).$$
Suppose the statement holds for $P_n(x)$, i.e.
$$P_n(x)=a_{n-1}x^{n-1}+\sum_{j=0}^{n-2} a_j x^j$$
where $a_{n-1}$ is a positive integer and the $a_j$'s ($j\le n-2$) are nonnegative integers. Then
\begin{equation}\label{eqn:coefficients}
P_{n+1}(x)=2 a_{n-1} x^{n} + \sum_{j=0}^{n-2} (n-j+1) a_j x^{j+1}+\sum_{j=0}^{n-1} j a_j x^{j-1}.
\end{equation}
Therefore $P_{n+1}(x)$ is a polynomial of degree $n$ whose coefficients are nonnegative integers. By induction, this completes the proof of the first part of the lemma.

The fact that $P_n(x)$ consists of either the even powers $1,x^2,\cdots,x^{n-1}$ or the odd powers powers $x,x^3,\cdots,x^{n-1}$ (depending on whether $n$ is odd or even) follows easily from the recursion formula \eqref{eqn:coefficients} and induction.
\end{proof}

We believe that the $N(\q)$'s the Proposition \ref{prop:convexity} should not be present, but we have not been able to remove them. By examining $h_\infty''(1/2)$ in its dependence on $\q$, we make the following conjecture, where $\zeta(s)$ denotes the Riemann zeta function.

\begin{conj}
The function
$$F(s)=2(s+1)(2^{s+2}-1)\zeta(s+2)-2\pi^2 (2^s-1)\zeta(s),\ 1<s<\infty.$$
is strictly increasing. In particular $s=2$ is the unique zero of $F(s)$.
\end{conj}

\begin{figure}[h]
\includegraphics[scale=0.4]{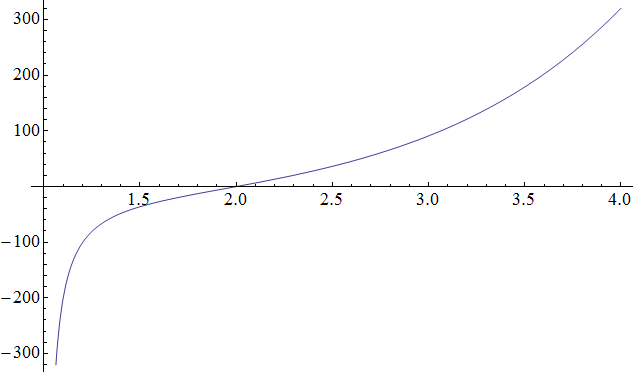}
\caption{A graph of $F(s)$ for $1<s<4$.}
\end{figure}

\newpage

\section{$L^p$ space}\label{sec:Lp}
Let $1\le p\le\infty$. We can also consider the transfer operator $L$ on the Lebesgue space $L^p(\mathbb T)$. Here we consider only the case $f(t)=|\cos(\pi t)|^\q$. Other cases can be treated similarly.

Let $d\ge 2$ be an integer. Let $L_\q=L: L^p(\mathbb T)\rightarrow L^p(\mathbb T)$ be the operator given by 
\begin{align*}
(L_\q u)(t)
=\frac{1}{d}\sum_{i=0}^{d-1}
\left |\cos\left (\pi\frac{t+i}{d}\right )\right |^\q
u\left (\frac{t+i}{d}\right ).
\end{align*}
Then $L_\q$ defines a bounded linear operator on $L^p(\mathbb T)$.
The adjoint of $L_\q$ is given by 
$T_\q=L_\q^*: L^{p'}(\mathbb T)\rightarrow L^{p'}(\mathbb T)$ (where $p'=p/(1-p)$),
\begin{align*}
(T_\q x)(t)
=|\cos(\pi t)|^\q x(d\cdot t).
\end{align*}
Notice that, for $p'<\infty$,
\begin{align*}
\|T_\q x\|_{p'}^{p'}
&=\int_{\mathbb T} |\cos(\pi t)|^{\q p'} |x(d\cdot t)|^{p'}dt\\
&=\int_{\mathbb T} (L_{\q p'}1)(t) |x(t)|^{p'}dt.
\end{align*}
Thus
\begin{align*}
\|L_\q\|_{p\rightarrow p}
&=\|T_\q\|_{p'\rightarrow p'}\\
&=\left (\sup_{\|x\|_{p'}\le 1} \|T_\q x\|_{p'}^{p'} \right )^{1/p'}\\
&=\left (\sup_{\||x|^{p'}\|_{1}\le 1}
\int_{\mathbb T} (L_{\q p'}1)(t) |x(t)|^{p'}dt \right )^{1/p'}\\
&=\|L_{\q p'}1\|_\infty^{1/p'}.
\end{align*}
By Lemma \ref{S:l1}, the function $h_1(t)=(L_{\q p'}1)(t)$ attains its maximum at either $t=0$ or $t=1/2$ depending on the value of $\q p'$. (Note that the function $h_1$ for $f=|\cos(\pi t)|^\q$ and that for  $f=|\sin (\pi t)|^\q$ differ only by a translation of $d/2$.) In particular, we obtain an explicit formula for the operator norm
$$\|L_\q\|_{p\rightarrow p}=\max\left \{(L_{\q p'}1)(0),(L_{\q p'}1)(1/2)\right \}.$$
More generally, for any $n\in\N$, the same argument as above gives
$$\|L_\q^{n}\|_{p\rightarrow p}=\|L_{\q p'}^{n}1\|_\infty^{1/p'}.$$
Therefore, to compute the spectral radius of $L_\q$ on $L^p(\mathbb T)$, it suffices to find
\begin{align*}
\rho_p(L_\q):=\lim_{n\rightarrow\infty} \|L_\q^n\|^{1/n}_{p\rightarrow p}
&=\lim_{n\rightarrow\infty}\|L_{\q p'}^{n}1\|_\infty^{1/(np')}\\
&=\lim_{n\rightarrow\infty}\|h_n\|_\infty^{1/(np')}
\end{align*}
where $h_n(t)=(L_{\q p'}^{n}1)(t)$. However, the last expression is the $p'$th root of the spectral radius of $L_{\q p'}$ on $C(\mathbb T)$. So we obtain the following.

\begin{prop} For $p>1$, we have
$$\rho_p (L_\q)=\big[\rho (L_{\q p'})\big]^{1/p'}.$$
\end{prop}

Similar as in Section \ref{sec:binary}, in the special case $d=2$, we can find eigenfunctions of $L_\q$ in $L^p(\mathbb T)$ explicitly. We consider two different cases.

Case 1: $\q p'\le 1$.  In this case we have, by Theorem \ref{2:t},
$$\rho (L_{\q p'})=2^{-\q p'},$$
and so
$$\rho_p (L_\q)=2^{-\q}.$$
Since $\q\le 1$, the spectral radius of $L_\q$ on $L^p(\mathbb T)$ coincides that on $C(\mathbb T)$. In particular, we have the same eigenfunction
$$u(t)=|\sin(\pi t)|^\q\in L^p(\T)$$
corresponding to the eigenvalue $\lambda=2^{-\q}$.

Case 2: $\q p'>1$.  In this case we have
$$\rho (L_{\q p'})=\frac12,$$
and so
$$\rho_p (L_\q)=2^{-1/p'}.$$
Note that $1/p'<\q$. Following the same idea as in Section \ref{sec:binary}, we consider functions of the form
$$u_s(t)=|\sin(\pi t)|^\q G(s,t)$$
where $s>1$ and $G(s,t)=\zeta(s,t)+\zeta(s,1-t)$ 
\footnote{More generally, one can take $G(s,t)$ to be linear combinations of $\zeta(s,t)$ and $\zeta(s,1-t)$.}
is as in \eqref{3:Gst}. Since
$$\zeta(s,t)\sim t^{-s},\quad \text{as }t\rightarrow 0^+,$$
we have that $u_s\in L^p(\mathbb T)$ if and only if $(s-\q)p<1$, i.e.
$$s<\q+\frac{1}{p}.$$
Since $\q p'>1$ exactly when $\q+\frac{1}{p}>1$,
we can take
$$s=\q+\frac{1}{p}-\varepsilon$$
for sufficiently small $\varepsilon>0$ to obtain an eigenfunction in $L^p(\mathbb T)$ corresponding to the eigenvalue
$$2^{-\q+(s-1)}=2^{-1/p'-\varepsilon}.$$
Therefore, as $\varepsilon\rightarrow 0$, $u_s(t)$ gives an `approximate' eigenfunction corresponding to $\rho_p (L_\q)=2^{-1/p'}$.
Note that when $\varepsilon=0$, $u_s(t)$ gives an eigenfunction in the Lorentz space $L^{p,\infty}(\mathbb T)$.
 
\section{An application to Fourier multipliers}\label{sec:application}
In this section, we present an application to some Bochner-Riesz type multipliers introduced by Mockenhaupt in \cite[Section~4.3]{Mockenhaupt}. Let $E\subset\R$ be the middle-third Cantor set obtained from dissecting the interval $[-1/2,1/2]$, and let $\mu$ be the Cantor measure on $E$. It is well known that
$$\dim E=\alpha:=\frac{\log 2}{\log 3}$$
and that the Fourier transform of $\mu$ is given by
\begin{align}\label{fourier-cantor}
\hat{\mu}(x)=\int_{\mathbb R} e^{-\pi i x \xi}d\mu(\xi)=\prod_{j=1}^\infty {\cos(\pi 3^{-j} x)}.
\end{align}
Let $\chi\in C_c^\infty(\R)$ be a bump function with $\hat{\chi}\ge 0$. For $\delta>0$, let
\begin{align*}
m_\delta(\xi)
=\frac{\chi(\cdot)}{|\cdot|^{\alpha-\delta}}*\mu(\xi)
=\int_\R \frac{\chi(\xi-\eta)}{|\xi-\eta|^{\alpha-\delta}}d\mu(\eta).
\end{align*}
Note that $m_\delta$ defines a bounded function only when $\delta>0$. In particular, $m_\delta$ is an $L^2$-Fourier multiplier if and only if $\delta>0$. 

\begin{thm}\label{multiplier}
$m_\delta$ is an $L^1$-Fourier multiplier if and only if
$$\delta>\frac{\log 2}{\log 3}+\frac{\log c(1)}{\log 3}=0.236\cdots$$
where $c(1)$ is as in Section \ref{sec:special-weights} (with $d=3$).
\end{thm}

\begin{proof}
Recall that an $L^p$-Fourier multiplier is a function $m(\xi)$ such that
\begin{equation}\label{eq:multiplier}
\|\mathcal{F}^{-1}\big(m(\xi)\hat f(\xi)\big)\|_{L^p(\R)}\le C \|f\|_{L^p(\R)}
\end{equation}
holds for a constant $C$ independent of $f$, where $\mathcal{F}^{-1}$ denotes the inverse Fourier transform. In the case $p=1$, this is equivalent to $\widehat {m}$ being a finite measure. If $\alpha-\delta\le 0$, it is easy to see that this is the case with $m=m_\delta$. If $\alpha-\delta>0$, then we have
$$\widehat{m}_\delta(x)=c\cdot\left (\hat{\chi}*{|\cdot|^{\alpha-\delta-1}}\right )(x)\cdot \hat{\mu}(x)$$
for some constant $c$. Thus, ${m}_\delta$ is an $L^1$-Fourier multiplier if and only if 
\begin{align*}
\int_{\R} |\widehat m_\delta(x)|dx
&=\int_{|x|\le 3} |\widehat m_\delta(x)|dx + \sum_{k=1}^\infty \int_{3^k<|x|\le 3^{k+1}} |\widehat m_\delta(x)|dx\\
&\approx 1 + \sum_{k=1}^\infty 3^{(\alpha-\delta-1)k} \int_{3^k}^{3^{k+1}} \prod_{j=1}^\infty |{\cos(\pi 3^{-j} x)}|dx\\
&<\infty
\end{align*}
where we have used 
$$\hat{\chi}*{|\cdot|^{\alpha-\delta-1}}(x)\approx |x|^{\alpha-\delta-1},\ \text{as } |x|\rightarrow\infty$$
and \eqref{fourier-cantor}. On the other hand, notice that
\begin{align*}
\int_{3^k}^{3^{k+1}} \prod_{j=1}^\infty |{\cos(\pi 3^{-j} x)}|dx
&=3^k \int_{1}^{3} \prod_{j=1}^\infty |{\cos(\pi 3^{k-j} x)}|dx\\
&=3^k \int_{1}^{3} |\hat \mu(x)| \prod_{j=0}^{k-1} |{\cos(\pi 3^{j} x)}| dx\\
&\approx 3^k \int_{0}^{1} \prod_{j=0}^{k-1} |{\cos(\pi 3^{j} x)}| dx
\end{align*}
where in the last line we have used periodicity and the fact that $|\hat\mu(x)|$ is bounded below on the interval $[2,3]$. Now by Theorem \ref{L:t4}(b), we know that
$$\int_{0}^{1} \prod_{j=0}^{k-1} |{\cos(\pi 3^{j} x)}| dx\approx c(1)^k.$$
Therefore 
$$\int_{\R} |\widehat m_\delta(x)|dx<\infty$$
if and only if
$$\sum_{k=1}^\infty 3^{(\alpha-\delta-1)k} 3^k c(1)^k<\infty,$$
which is equivalent to
$$\delta>\frac{\log 2}{\log 3}+\frac{\log c(1)}{\log 3}.$$
This completes the proof.
\end{proof}

Since $m_\delta$ is compactly supported, we can choose $f\in L^p(\R)$ in \eqref{eq:multiplier} such that $\hat f\equiv 1$ on the support of $m=m_\delta$, and get $\widehat m_\delta\in L^p(\R)$ as a necessary condition for $m_\delta$ to be an $L^p$-Fourier multiplier. By the same argument as above, this leads us to
$$\delta>\delta(p):=\frac{\log 2}{\log 3}-1+\frac{1}{p}+\frac{\log \big(c(p)^{1/p}\big)}{\log 3}.$$

\begin{figure}[h]
\includegraphics[scale=0.4]{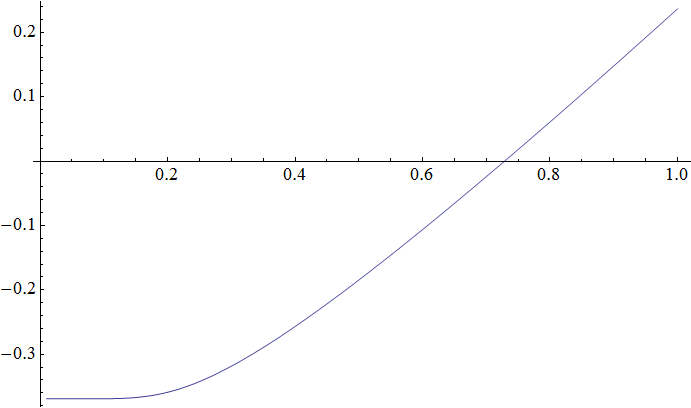}
\caption{A graph of $\delta(p)$ as a function of $1/p\in (0,1)$.}
\end{figure}

\bibliographystyle{abbrv} 
\bibliography{bibliography}
\nocite{JanardhanRosenblumStrichartz1992}

\end{document}